\documentclass[12pt]{amsart} 

\usepackage{url}           
\usepackage{amscd,amssymb,mathtools,enumerate}
\usepackage[plainpages,urlcolor=blue,backref]{hyperref}

\theoremstyle{plain}
\newtheorem{thm}{Theorem}
\newtheorem{cor}{Corollary}
\newtheorem{lem}{Lemma}

\newcommand{\A}{{\mathcal A}}
\newcommand{\Ann}{\mathrm{Ann}}
\newcommand{\f}{{\bf f}}
\newcommand{\C}{\mathbb{C}}
\newcommand{\p}{\partial}
\newcommand{\PP}{\mathbb{P}}
\newcommand{\V}{{\mathcal V}}

\begin{document}
\date{}

  \title[Strong Lefschetz Property and Associated Forms]{On Strong Lefschetz Property of 0-dimensional Complete Intersections and Associated Forms}

   \author[ZHENJIAN WANG]{ZHENJIAN WANG}
   \address{Hefei National Laboratory, Hefei 230088, China}
   \email{wzhj@ustc.edu.cn}



\begin{abstract} 
We prove that a homogeneous 0-dimensional complete intersection satisfies the Strong Lefschetz Property (SLP) in degree 1 if and only if its associated form has nonzero Hessian. This result is essentially known in the literature, but our proof differs from previous ones. 
We then investigate properties of the associated forms and show that they can always be reconstructed from their partial derivatives, in particular from their Jacobian ideals. 
As an application of our method, we prove that every smooth homogeneous polynomial which is not of Sebastiani-Thom type is uniquely determined by any nontrivial graded component of its Jacobian ideal.
\end{abstract}

\maketitle


\section{Statements of main results}

In \cite{DIN25}, the authors establish that a homogeneous 0-dimensional complete intersection satisfies the Strong Lefschetz Property (SLP) in degree 1 provided that its associated form has a nonzero discriminant. Here, we show that this condition can be weakened: the necessary and sufficient condition is that the Hessian of the associated form is nonzero.

Following the notation of \cite{DIN25}, let $S=\C[x_1, \ldots,x_n]$ be the polynomial ring with the usual grading and $n \geq 3$, and $f_j \in S_{d_j}$ be a homogeneous polynomial of degree $d_j \geq 2$ for $j=1, \ldots,n$ such that the ideal
\[
J(\f)=(f_1, \ldots ,f_n)
\]
is a 0-dimensional complete intersection.
The quotient $M(\f)=S/J(\f)$ is a graded Artinian Gorenstein algebra with socle degree
\begin{equation}
\label{eq0}
T=\sum_{j=1}^nd_j-n.
\end{equation}

To $\f$, one can associate a homogeneous polynomial $\mathrm{A}_\f$, called \textit{associated form} of $\f$, in the polynomial ring $R:=\C[y_1,\ldots,y_n]$, where $y_i$ are dual variables to $x_i$~\cite{AI14,AI18,EI13}. More precisely, denote by $\omega\colon \mathrm{Soc}(M(\f))\to  \C$ the linear isomorphism given by the condition $\omega(\overline{\mathrm{Jac}(\f)})=1$, where $\mathrm{Jac}(\f)$ is the determinant of the Jacobian matrix of the induced polynomial mapping $\f: \C^n \to \C^n$. The associated form $\mathrm{A}_\f$ is defined by the formula 
\begin{equation}
\label{eqAS1} 
\mathrm{A}_\f(y_1,\ldots,y_n)=\omega((y_1\overline{x_1}+\dots+y_n\overline{x_n} )^T),
\end{equation}
where $\overline{x_i}\in M(\f)$ is the image of $x_i$.

For a homogeneous polynomial $F\in R$ of degree $T$, its \textit{apolar ideal} of $F$ is defined by 
\[ 
I_F=\mathrm{Ann}(F)=\left\{ g\in  S \ : \ g\circ F=0\right\}\subseteq S,
\]
where the action of $S$ on $R$ is obtained by identifying $x_i$ with
$\partial /\partial y_i$.
Moreover, for any finitely generated graded Artinian Gorenstein $\C$-algebra $S/J$ with socle degree $T$, there exists a homogeneous polynomial $F\in R_T$, unique up to scaling, such that $J=I_F$; see~\cite[Lemma 2.12]{IK99}. Any such $F$
is called a \textit{Macaulay inverse system} for $S/J$. It is known that 
the associated form $\mathrm{A}_\f(y_1,\ldots,y_n)\in R_T$ is a Macaulay inverse system for the Artinian algebra $M(\f)$, see~\cite[Proposition 3.2]{AI14}.

In this note, we provide a necessary and sufficient condition on the ideal $J(\f)$ for the Artinian algebra $M(\f)$ to satisfy SLP in degree 1, expressed in terms of the Hessian of its associated form.

\begin{thm}\label{thm1}
Let $M(\f)$ be a 0-dimensional homogeneous complete intersection.
Then the strong Lefschetz property holds in degree 1 for $M(\f)$; that is, for a generic linear form $\ell \in S_1$, the multiplication map
\[
\ell^{T-2}:M(\f)_1 \to M(\f)_{T-1}
\]
is an isomorphism, if and only if its associated form $\mathrm{A}_\f$ has nonzero Hessian.
\end{thm}

As in \cite[Remark 1.2]{DIN25}, let $f\in S_d$ be a reduced homogeneous polynomial of degree $d\geq 3$. By Euler's formula, the singular subscheme of the reduced hypersurface $V(f)\subseteq \PP^{n-1}$ is defined by the gradient ideal 
\[ 
J(f):=\left(f_1=\frac{\partial f}{\partial x_1},\dots,f_n=\frac{\partial f}{\partial x_n}\right).
\]
The graded $\C$-algebra $M(f):=S/J(f)$ is called the Milnor algebra of $f$. If $f$ defines a smooth hypersurface in $\PP^{n-1}$, then $J(f)$ is generated by a regular sequence; thus, $M(f)$ is an Artinian Gorenstein graded $\C$-algebra with socle degree $T=n(d-2)$. We denote by
\[
A_f(y_1,\ldots,y_n)\in R_T
\]
the  Macaulay inverse system for the Milnor algebra $M(f)$.
This setting is a special case of 0-dimensional homogeneous complete intersections described above.
 Hence Theorem \ref{thm1} immediately yields the following.

\begin{cor}\label{cor1}
Let $f$ be a polynomial of degree $d\geq 3$ defining a smooth projective hypersurface. Then 
the strong Lefschetz property holds in degree 1 for the Milnor algebra $M(f)$ if and only if the associated form $A_f$ has nonzero Hessian.  
\end{cor}

We were made aware that Theorem \ref{thm1} has appeared in prior works (see e.g. \cite[Theorem 3.1]{MW09},\cite[Corollary 2.14 (2)]{DGI20}). Therefore, we present a self-contained proof of this known result using a different approach.

The second part of this article concerns properties of the associated forms $\mathrm{A}_\f$.
To present the results, we first introduce the following notation. For a homogeneous polynomial $G\in R_T$ and an integer $0\leq k\leq T-1$, denote by $E_k(G)$ the vector subspace of $R_{T-k}$ spanned by all partial derivatives of $G$ of order $k$ (cf. \cite{Wan20}).

\begin{thm}\label{thm2}
Let $M(\f)$ be a 0-dimensional homogeneous complete intersection. For any $0\leq k\leq T-\max\limits_{1\leq j\leq n}d_j$, the associated form $\mathrm{A}_\f$ can be reconstructed from $E_k(\mathrm{A}_\f)$; that is, if $G\in R_T$ satisfies $E_k(G)=E_k(\mathrm{A}_\f)$, then $G\in\C^* \mathrm{A}_\f$.
\end{thm}

In particular, $\mathrm{A}_\f$ is determined by $J(\mathrm{A}_\f)_{T-1}$. Hence, $\mathrm{A}_\f$ can always be reconstructed from its Jacobian ideal (cf. \cite{Wan15}).

Theorem \ref{thm2} immediately yields the following.

\begin{cor}\label{cor2}
Let $f\in S$ be a polynomial of degree $d\geq 3$ defining a smooth projective hypersurface. Then for any $0\leq k\leq (n-1)(d-2)+1$, the associated form $A_f$ can be reconstructed from  $E_k(A_f)$.
\end{cor}

Finally, as an application of our method, we prove the following result, which can be seen as a generalization of the results in \cite{UY09, Wan15}.

\begin{cor}\label{cor3}
Let $f\in S$ be a polynomial of degree $d\geq 3$ defining a smooth projective hypersurface, and let $d-1\leq k\leq T=n(d-2)$. Assume that $f$ is not of Sebastiani-Thom type and moreover that $J(f)_k=J(g)_k$ for some $g\in S_d$. Then $g\in\C^*f$. In other words, $f$ can be reconstructed from the graded component $J(f)_k$ of the Jacobian ideal of $f$.
\end{cor}

\section{Proof of Theorem \ref{thm1}}
Consider the Veronese morphism
\[
v: \PP(S_1)\to\PP(S_{T-1}).
\]
Since the Veronese variety
\[
\V=v(\PP(S_1))
\]
 spans the whole projective space $\PP(S_{T-1})$ and $\PP(J(\f)_{T-1})$ is a proper linear subspace of $\PP(S_{T-1})$, the set
\[
\A=v^{-1}\left(\V\cap\PP(J(\f)_{T-1})\right)
\]
is a proper algebraic subset of $\PP(S_1)$. 

Let
\[
\pi:  \PP(S_{T-1}) \setminus \PP(J(\f)_{T-1}) \to \PP^{n-1}
\]
be the linear projection with center $\PP(J(\f)_{T-1})$; it is a regular morphism.
Now consider the composition
\[
\phi=\pi \circ v: \PP(S_1)\setminus\A \to \PP^{n-1}
\]
which is also regular.  

As indicated in the proof of Theorem~1.5 of \cite{DIN25}, the algebra $M(\f)$ has SLP in degree 1 exactly when $\phi$ is an immersion at some point $\ell\in \PP(S_1)\setminus\A$. Indeed, take any linear form $\ell\in\PP(S_1)\setminus\A$. The projective tangent space to the Veronese variety $\V$ at the point $p=\ell^{T-1}$ is given by
\begin{equation}
\label{eq3} 
T_p\V= \left\{\ell^{T-2}\ell_1 \  : \ \ell_1 \in S_1\right\} \subset \PP(S_{T-1}),
\end{equation}
see for instance \cite[Section 1]{CGG02}. Consequently, the kernel of the differential $d\phi_\ell$ 
consists of those $\ell_1\in S_1\setminus\C\ell$ for which $\ell^{T-2}\ell_1$ belongs to the span of $\ell^{T-1}$ and $J(\f)_{T-1}$. 

If SLP in degree 1 fails for $M(\f)$, then for every linear form $\ell\in\PP(S_1)\setminus\A$, there exists a nonzero $\ell_1\in S_1$ such that $\ell^{T-2}\ell_1\in J(\f)_{T-1}$. Note that such $\ell_1$ cannot lie in $\C\ell$ because $\ell^{T-1}\notin J(\f)_{T-1}$. Hence $\phi$ is nowhere an immersion. 

Conversely, assume that SLP holds in degree 1 but $\phi$ is not an immersion at any point. Then for a generic linear form $\ell$, we can find $\ell_1\in S_1\setminus\C\ell$ and a scalar $a\in\C$ with $(\ell_1-a\ell)\ell^{T-2}\in J(\f)_{T-1}$. Then SLP forces $\ell_1=a\ell$, a contradiction.

The morphism $\phi$ can be written explicitly as
\[
\phi=[F_1,\ldots, F_n]
\]
where $F_1,\ldots, F_n\in R_{T-1}$ form a basis of the annihilator $\mathrm{Ann}(J(\f)_{T-1})$, defined by
\[
\mathrm{Ann}(J(\f)_{T-1})=\{F\in R_{T-1}\ :\ g\circ F=0\text{ for all }g\in J(\f)_{T-1}\}.
\]
In fact, identifying $S$ with $R$ via $x_i\longleftrightarrow y_i$, the polar pairing induces a Hermitian inner product on $S_{T-1}$, given by
\[
\langle g, F\rangle=g\circ\overline{F},
\]
where $\overline{F}$ denotes the complex conjugation of $F$; that is,
\[
\overline{F}=\sum \overline{a_\alpha} x^\alpha\quad\text{if}\quad F=\sum a_\alpha x^\alpha.
\]
Choose a basis $G_1,\ldots, G_n$ of $\mathrm{Ann}(J(\f)_{T-1})$ satisfying the orthogonality condition $G_i\circ\overline{G_j}=\delta_{ij}$. Observe that $\overline{G_1},\ldots,\overline{G_n}$ span the orthogonal complement of $J(\f)_{T-1}$ with respect to this Hermitian inner product. Consequently, the morphism $\phi$ can be expressed as
\[
\phi(a)=[g_1(a),\ldots, g_n(a)],\ a=[a_1,\ldots, a_n]
\]
where $g_i(a)$ are determined by
\[
(a_1x_1+\dots+a_nx_n)^{T-1}\equiv g_1(a)\overline{G_1}+\dots+g_n(a)\overline{G_n}\mod J(\f)_{T-1}.
\]
Applying the polar pairing with $G_i$ to both sides and using the fact that $G_i\circ J(\f)_{T-1}=0$, we obtain
\[
g_i(a)=G_i\circ(a_1x_1+\dots+a_nx_n)^{T-1}=(T-1)!G_i(a).
\]

Thus it suffices to show that the basis $F_1,\ldots, F_n$ can be chosen algebraically independent over $\C$. 

Recall that $J(\f)=\mathrm{Ann}(\mathrm{A}_\f)$. By~\cite[Lemma 2.15]{IK99}, we have $J(\f)_{T-1}=\mathrm{Ann}(\mathrm{A}_\f)_T:S_1$, where
\[
\mathrm{Ann}(\mathrm{A}_\f)_T:S_1=\left\{g\in S_{T-1}\ :\  gS_1\subseteq\mathrm{Ann}(\mathrm{A}_\f)_T\right\}.
\]
 Since
\[
(x_i g)\circ \mathrm{A}_\f=0\Longleftrightarrow g\circ\left(\frac{\partial \mathrm{A}_\f}{\partial y_i}\right)=0,\ 1\leq i\leq n,
\]
it follows that $\mathrm{Ann}(J(\f)_{T-1})$ is spanned by the first-order partial derivatives of the associated form.  These derivatives are algebraically independent if and only if the Hessian of $\mathrm{A}_\f$ is nonzero. The proof is now complete.

\section{Proof of Theorem \ref{thm2}}

Recall that the polar action of $S=\C[x_1,\ldots, x_n]$ on $R=\C[y_1,\ldots, y_n]$ is given by
\[
S\times R\to R,\qquad (g,F)\mapsto g\circ F:=g\left(\frac{\p}{\p y_1},\ldots,\frac{\p}{\p y_n}\right)F.
\]
Given a vector subspace $E\subset R_k$, denote by
\[
E^\perp=\left\{g\in S_k\  :\  g\circ F=0\text{ for all } F\in E\right\}.
\]
In addition, for two vector spaces $E_1\subset S_i$ and $E_2\subset S_j$, set
\[
E_1:E_2=\left\{g\in S_{i-j}\  :\  gg_2\in E_1\text{ for all } g_2\in E_2\right\},
\]
and
\[
E_1\cdot E_2=\left\{g_1g_2\in S_{i+j}\ :\  g_1\in E_1,\  g_2\in E_2\right\}.
\]

Given any $F\in R_T$, from the identity
\[
g\circ\left(\frac{\p^\alpha F}{\p y^\alpha} \right)=\left(x^\alpha g\right)\circ F
\]
we obtain that
\[
E_k(F)^\perp=\left(\C F\right)^\perp:S_k=\Ann(F)_{T}:S_k.
\]
By~\cite[Lemma 2.15]{IK99}, we have $\Ann(F)_T:S_k=\Ann(F)_{T-k}$. Therefore, 
\begin{equation}
\label{eq:perp1}
E_k(F)^\perp=\Ann(F)_{T-k}.
\end{equation}

\begin{lem}
Let $M(\f)$ be a 0-dimensional homogeneous complete intersection and $\mathrm{A}_\f$ its associated form. Then $\mathrm{A}_\f$ is not a cone; that is, the first partial derivatives $\frac{\p \mathrm{A}_\f}{\p y_i},i=1,\ldots, n$ are linearly independent over $\C$.
\end{lem}
\begin{proof}
It suffices to show that $\dim E_1(\mathrm{A}_\f)=n$. Indeed, we have
\[
\dim E_1(\mathrm{A}_\f)=\mathrm{codim}\, E_1(\mathrm{A}_\f)^\perp
\]
and, by \eqref{eq:perp1}, $E_1(\mathrm{A}_\f)^\perp=\Ann(\mathrm{A}_\f)_{T-1}$. By the definition of the associated form $\mathrm{A}_\f$, we have $\Ann(\mathrm{A}_\f)=J(\f)$. Hence, 
\[
\dim E_1(\mathrm{A}_\f)=\mathrm{codim}\, J(\f)_{T-1}=\dim M(\f)_{T-1}=\dim M(\f)_1.
\]
Since $d_j\geq2$ for $j=1,\ldots, n$, it follows that $\dim M(\f)_1=n$.
\end{proof}

\subsection{Proof of Theorem \ref{thm2}}
Let $M(\f)$ be a 0-dimensional homogeneous complete intersection and $\mathrm{A}_\f$ its associated form. Assume $1\leq k\leq T-\max\limits_{1\leq j\leq n}d_j$. Let $G\in R_T$ satisfy $E_k(G)=E_k(\mathrm{A}_\f)$.

From \eqref{eq:perp1}, we obtain
\[
\Ann(G)_{T-k}=E_k(G)^\perp=E_k(\mathrm{A}_\f)^\perp=\Ann(\mathrm{A}_\f)_{T-k}=J(\f)_{T-k}.
\]
Since $k\leq T-\max\limits_{1\leq j\leq n}d_j$, we have $J(\f)_T=S_k\cdot J(\f)_{T-k}$. It follows that
\[
J(\f)_T=S_k\cdot \Ann(G)_{T-k}\subseteq\Ann(G)_T.
\]
On the other hand, both $J(\f)_T$ and $\Ann(G)_T$ have codimension 1 in $S_T$. Hence they must be equal: $J(\f)_T=\Ann(G)_T$.

Note that $(\C \mathrm{A}_\f)^\perp=J(\f)_T$ and $(\C G)^\perp=\Ann(G)_T$. Consequently, we deduce that $\C \mathrm{A}_\f=\C G$, and thus $G\in\C^*\mathrm{A}_\f$. The proof of Theorem \ref{thm2} is now complete.

\section{Application: Proof of Corollary \ref{cor3}}

Let $f\in S$ be a polynomial of degree $d\geq 3$ defining a smooth projective hypersurface. Let $d-1\leq k\leq T=n(d-2)$ and let $g\in S_d$ satisfy $J(f)_k=J(g)_k$.

First of all, $g$ also defines a smooth projective hypersurface. Indeed, the singular locus of the hypersurface $V(g)\subseteq\PP^{n-1}$ is the zero set of $J(g)_k=J(f)_k$, which is empty by the smoothness of $f$.

In particular, the associated form $A_g$ for $M(g)$ is well-defined.

Now by \eqref{eq:perp1}, the condition $J(f)_k=J(g)_k$ is equivalent to the equality $E_{n(d-2)-k}(A_f)^\perp=E_{n(d-2)-k}(A_g)^\perp$. Thus $E_{n(d-2)-k}(A_f)=E_{n(d-2)-k}(A_g)$. It follows from Corollary \ref{cor2} that $A_f\in\C^* A_g$. Consequently,
\[
J(f)=\Ann(A_f)=\Ann(A_g)=J(g).
\]
Then the conclusion of Corollary \ref{cor3} follows from \cite[Lemma 3]{UY09} or \cite[Corollary 1.3]{Wan15}.

\section*{Acknowledgements}

We would like to thank Professor A. Dimca for useful discussions and suggestions. We also thank Hefei National Laboratory for their wonderful working atmosphere.


\end{document}